\documentclass[a4paper,10pt,oneside,onecolumn,preprint]{elsarticle}

\usepackage{graphicx}
\usepackage{amsmath}
\usepackage{amssymb}
\usepackage{subcaption}
\usepackage{tikz}
\usetikzlibrary{backgrounds}
\usetikzlibrary{calc,trees,positioning,arrows,fit,shapes}
\usepackage{float}
\usepackage{hyperref}
\usepackage{verbatim}
\usepackage{mathtools}

\newtheorem{theorem}{Theorem}[section]
\newtheorem{proposition}{Proposition}[section]
\newtheorem{corollary}{Corollary}[section]
\newtheorem{definition}{Definition}[section]
\newproof{proof}{Proof}
\newtheorem{lemma}{Lemma}[section]
\newtheorem{claim}{Claim}[section]

\newcommand{\di}{\Diamond}
\newcommand{\subf}[1]{\mathsf{sub}(#1)}
\newcommand{\Power}[1]{{  \mathcal{P} }(#1)}

\newcommand{\checkpath}{\mathtt{CheckPath}}
\newcommand{\PSPACE}{\textsc{PSpace}}

\newcommand{\SAW}{\mathfrak {C}}

\newcommand\pto{\mathrel{\ooalign{\hfil$\mapstochar$\hfil\cr$\to$\cr}}}

\title{Modal logic of planar polygons}
\author[rmi]{David Gabelaia \footnote[1]{Supported by Shota Rustaveli National Science Foundation grant \#DI-2016-25}}
\author[illc]{Kristina Gogoladze}
\author[rmi]{Mamuka Jibladze \footnotemark[1]}
\author[tsu]{\\Evgeny Kuznetsov \footnotemark[1]}
\author[ilps]{Maarten Marx}

\address[ilps]{ILPS, Informatics Institute, Universiteit van Amsterdam}
\address[illc]{ILLC, Institute for Logic Language and Computation ,Universiteit van Amsterdam}
\address[tsu]{TSU, Ivane Javakhishvili Tbilisi State University}
\address[rmi]{TSU, Andria Razmadze Mathematical Institute, Ivane Javakhishvili Tbilisi State University}

\begin{document}
\begin{frontmatter}

\begin{abstract}
We study the modal logic of the closure algebra $P_2$, generated by the set of all polygons in the Euclidean plane $\mathbb{R}^2$. We show that this logic is finitely axiomatizable, is complete with respect to the class of frames we call ``crown'' frames, is not first order definable, does not have the Craig interpolation property, and its validity problem is $\PSPACE$-complete.
\end{abstract}
\begin{keyword}
Modal Logic, topological semantics, planar polygons.
\end{keyword}
\end{frontmatter}

\section{Introduction}\label{intro}
Connections between modal logic and general topology were discovered and investigated in the works of McKinsey and Tarski in the 1940s (see, e.g., their seminal paper \cite{mck:mckinsey44}). Perhaps this timing, coupled with the arrival of the rival Kripke-style semantics in the 1950s has led to the slow progress in this area of research. Already from the beginning it was established that if the modal diamond is interpreted as the closure operator over a topological space, then the minimal modal logic is $\mathbf{S4}$. Moreover, for arbitrary $n\in \mathbb{N}$ the modal logic of the Euclidean space $\mathbb{R}^n$ is also $\mathbf{S4}$. This result can be interpreted to mean that the modal language is not expressive enough to distinguish the Euclidean spaces from each other; in particular the modal language is insensitive to dimension. One can increase the expressive power by extending the language (e.g. by adding the global modality, or the difference operator) or change the interpretation of the modal diamond to another topological operator (for instance - the limit operator, as suggested already by McKinsey and Tarski in \cite{mck:mckinsey44}). 
Another road to take is to restrict the \emph{valuations}, allowing the propositional letters to denote some \emph{well-behaved} subsets instead of arbitrary subsets of the space in question, where `well-behaved' might carry a topological meaning like `regular closed' or a geometrical meaning like `convex' or `polygonal'. In \cite{bezh:benthbezhgerke03, aiello:BB03} authors considered the collections of serial, convex, hyper-rectangular and chequered subsets of various Euclidean spaces and calculated the arising modal logics. In \cite{wolt:spat00} the regular closed regions of (Euclidean) spaces were considered from the mereotopological point of view and in \cite{pratt:pratt10} the polygonal regions were considered again from the point of view of mereotopology (see also \cite{zakh:conn13,zakh:rcc14} for the recent developments in this direction). In the latter line of work however, the boolean operations on regions do not always coincide with the usual set-theoretic operations. In particular, the meet of two regular closed regions may not coincide with their intersection. To phrase it algebraically, the propositional part of the modal language is being interpreted in the boolean algebra of regular closed subsets, rather than in the powerset boolean algebra. This differs from our approach --- we work with the powerset algebra, as, for example, in \cite{bezh:benthbezhgerke03}.

In this paper, we consider pure modal language interpreted over the two-dimensional Euclidean plane $\mathbb{R}^2$ in such a way that the propositional letters denote only the so-called \emph{polygonal} regions of the plane.

We introduce the structure $\mathfrak P_2=(\mathbb{R}^2, P_2)$ which we call the \emph{polygonal plane}.
Here $P_2$ is the boolean subalgebra of the powerset of $\mathbb{R}^2$ generated by all half-planes; we call elements of $P_2$ \emph{polygons}. The algebra $P_2$ turns out to be a closure subalgebra of the full closure algebra of all subsets of the real plane.

We interpret the modal language over the polygonal plane using valuations $\nu:\textsc{Prop}\to P_2$ assigning to variables polygons from $P_2$. These are then extended to all modal formulas in a standard way.

Our main object of study is the modal logic $\mathbf{PL_2}$ given by all modal formulas valid over $\mathfrak P_2$.

\

\emph{Organisation of the paper}

In section \ref{preliminaries}, we give main definitions and some basic notions from the area of topological and algebraic logic, which are related to the main matter and are employed in proofs. As the end of the  preliminary section \ref{preliminaries} we define the notion of first-order definability of class of models and give known in a literature tools to investigate the question of firs-order definability.  Thereafter in section \ref{fmp}, we prove the finite model property for $\mathbf{PL_2}$. In section \ref{axiomatization}, we show that the logic is finitely axiomatizable and hence, decidable. In section \ref{complexity}, the complexity of the satisfiability problem for $\mathbf{PL_2}$ is shown to be in $\textsc{PSpace}$-complete. In section \ref{interpolation} is pointed out and shown by constructing particular counterexample, that the logic of planar polygons does not have Craig interpolation property. In section \ref{frstrddfnblt}, we show that the logic $\mathbf{PL_2}$ is not first order definable.

\section{Preliminaries}\label{preliminaries}

Syntactically we are dealing with the basic modal language with the set of countably many propositional letters $\textsc{Prop}$ and the formulas built in the usual way using the propositional connectives together with unary modal connectives $\Diamond$ and $\Box$.

$$\varphi := p~|~ \bot ~|~ \neg \phi ~|~ \phi_{1} \vee \phi_{2} ~|~ \Diamond \phi ~$$
where $p$ ranges over elements of $\textsc{Prop}$, and $\phi$, $\phi_{1}$, $\phi_{2}$ are formulas.
The additional connectives such as $\wedge$, $\to$, $\leftrightarrow$ and the modality $\Box$ are defined as usual.

\ 

Semantically our object of study is the \emph{polygonal plane} mentioned above. In more detail, we define it as follows. Consider the regions of the plane obtained by the intersections of finitely many half-planes and generate the boolean algebra using the set-theoretic operations from these regions. An arbitrary member of the obtained boolean algebra is called a \emph{polygon} and the collection of all polygons is denoted by $P_2$. The structure $\mathfrak P_2=(\mathbb{R}^2, P_2)$ is called the \emph{polygonal plane}.

A typical bounded member of $P_2$ is a finite union of (open) $n$-gons, line segments and points. In other words, we consider as entities not only the 2-dimensional $n$-gons, but also their boundaries, i.e. `polygons' of lower dimension.

The algebra $P_2$ is a closure subalgebra of the full closure algebra of all subsets of the real plane (see \ref{t:cloclo}).

To interpret the modal language over the polygonal plane, we allow for valuations $\nu:\textsc{Prop}\to P_2$ to range over polygons only. The valuations are extended to arbitrary modal formulas using the set-theoretic counterparts for the propositional connectives, interpreting $\Diamond$ as the topological closure, and $\Box$ as topological interior operators.

The set of all valid modal formulas over $\mathfrak P_2$ is denoted by $\mathbf{PL_2}$.

\ 

Hereafter, we give some definitions of terms from topological and algebraic  logic ,  related and used in the paper. Experienced reader can skip this part with ease and without loss of clarity.
\subsection{Kripke frames, topological spaces and their morphisms}

A nonempty set $X$ together with a binary relation $R\subseteq X \times X $ is said to be a \emph{Kripke frame} and will be  denoted by $(X,R)$. To indicate that $(x,y)\in R$ holds we often write $xRy$ ($x$ sees $y$ by $R$); in such a case an element $y$ is called a successor of $x$, and $x$ - a predecessor of $y$ respectively. A subset $U\subseteq X$ is called \emph{upwards closed subset} (or \emph{up-set} for short) if it contains all successors of all its elements. A subset $F\subseteq X$ is called \emph{Downwards closed subset} (or \emph{down-set} for short) if it contains all predecessors of all its elements. It is easy to see that a complement of an up-set is a down-set and vice versa.  If $R$ is a relation on $X$, and $A\subseteq X$, the set $\{y\in X: \exists x \in A ~~s.t.~~ xRy\}$ actually is the set of successors of elements of $A$ and is denoted by $R(A)$; the set $\{y\in X: \exists x \in A ~~s.t.~~ yRx\}$ actually is the set of predecessors of elements of $A$ and is  denoted by $R^{-1}(A)$.

A  nonempty set $X$ together with a collection $\tau \subseteq \mathcal{P}(X)$ of subsets that is closed under finite intersections and arbitrary unions is said to be a \emph{topological space} and is denoted $(X,\tau)$. Often, when there is no ambiguity,  instead of  $(X,\tau)$ we write just $X$. The members of $\tau$ are called \emph{open} subsets, or simply opens. Their complements are called \emph{closed} subsets. (see \cite{enge:engelking77})

For a subset $A\subseteq X$, there exists the greatest open subset contained in $A$ (i.e. the union of all the opens contained in $A$) which is denoted by $\mathbb{I}A$ (read: ‘interior $A$’). Thus $\mathbb{I}$ is an operator over the subsets of the space $X$.

The closure operator, which is a dual of the interior operator, is defined by $\mathbb{C}A=-\mathbb{I}(-A)$ where $"-"$ stands for the set-theoretic complementation. Observe that $\mathbb{C}A$ is the
least closed set containing A (i.e. the intersection of all the closed sets including $A$ as subset). Thus $\mathbb{C}$ as well is an operator over the subsets of the space $X$ .

It is well known, that there is a close relationship between Kripke frames and topological spaces.  For given topological space $(X, \tau)$ one can define a relation $R_{\tau}$ on $X$ (\emph{specialization relation} corresponding to $\tau$) in the following way:  $xR_{\tau} y$ iff $x\in \mathbb{C}\{y\}$. It is easy to see, that defined relation is a \emph{reflexive} ($\forall x\in ~X (xR_{\tau}x)$) and \emph{transitive} ($\forall x,y,z \in ~X ((xR_{\tau}y, ~yR_{\tau}z)\Rightarrow xR_{\tau}z)$) relation on $X$.

On the other hand, for a given reflexive and transitive Kripke frame $(X,R)$ one can generate a topology $\tau_R$ on $W$ by declaring the up-sets of $W$ to be open subsets. It is easy to see that the topology thus generated will be of a special kind, in particular such that \emph{arbitrary intersections} of open sets are open as well (i.~e. for every point there is the least open subset containing it). Such spaces are called \emph{Alexandroff} topological spaces and indeed there is a one-to-one correspondence between Alexandroff spaces and reflexive-transitive Kripke frames; namely $R_{\tau_{R}}=R$ and  $\tau_{R_{\tau}}=\tau$ iff $R$ is a reflexive and transitive relation and $\tau$ is an Alexandroff topology.

Since the Kripke frames and topological spaces are already defined, we have to define the notion of their morphisms.

A map $f:(X_1 , R_1) \to (X_2 , R_2 )$ between Kripke frames is said to be \emph{monotone} if whenever $xR_1 y$, then $f(x) R_2 f(y)$.

A map $f:(X_1 , R_1) \to (X_2 , R_2 )$ between Kripke frames is said to be a \emph{p-morhism} (\emph{pseudo epimorphism}) if it is monotone, and whenever $f(x)R_2 z$, there exists $y\in X_1$ such that $x R_1 y$ and $f(y)=z$.

\begin{center}
\begin{tikzpicture}[>=latex]
\node (1) at (0,0) {\(y\)};
\node (2d) at (0,-1.5) {\(x\)};
\node (2r) at (1.5,0) {\(z\)};
\node (3) at (1.5,-1.5) {\(f(x)\)};

\draw[dashed,|->] (1) to node[above] {} (2r);
\draw[->] (2d) to node[left] {$R_{1}$} (1);
\draw[->] (3) to node[right] {$R_{2}$} (2r);
\draw[|->] (2d) to node[right] {} (3);
\begin{scope}[shift=($(3)!.3!(1)$)]
\end{scope}
\end{tikzpicture}
\end{center}

In case of topological spaces, instead of preserving relation we need to preserve topology. That is a map $f:X_1 \to X_2$ between topological spaces $X_1=(X_1,\tau_1)$ and  $X_2=(X_2,\tau_2)$ is said to be \emph{continuous} if whenever $O_2$ is open in $X_2$, i.e $O_2\in \tau_2$, then the set $f^{-1}(O_2)$ is open in $X_1$ i.e. $f^{-1}(O_2)\in\tau_1$. It worth to mention the well known fact, a map $f:(X_{1},R_{1})\to (X_{2},R_{2})$ between quasy-ordered sets is monotone if and only if it is continuous with respect to Alexandroff topology of up-sets.

A map $f:X_1 \to X_2$ between topological spaces $X_1=(X_1,\tau_1)$ and  $X_2=(X_2,\tau_2)$ is said to be \emph{open} if it sends open sets to open sets, i.e.:
 \begin{center}
 $O_1 \in \tau_1$ implies $f(O_1) \in \tau_2$.\\
 \end{center}
A map is called an \emph{interior map}, if it is both open and continuous. In case the map $f: X_1\to X_2$ is interior, its image $f(X_1)\subseteq X_2 $ is called interior image of $X_1$. A map $f:(X_1 , R_1) \to (X_2 , R_2 )$ between quasi-ordered sets is $p$-morphism if and only if   $f:(X_1 , \tau_{R_1}) \to (X_2, \tau_{R_2} )$ is interior map.

Interior maps play a special role in the topological semantics of modal logic, as well as $p$-mophisms in Kripke semantics, as we proceed to briefly recall in the following subsection.

\subsection{Topological semantics of modal logic}

Now we recall the topological semantics for the basic modal language. For reader interested in more extensive reference we suggest \cite{spatlog:AielloBenthem07}.

Nowadays, one of the best-known semantics for $\mathcal{ML}$ is the Kripke semantics. However, in this paper we study the topological semantics, according to which modal formulas denote regions in a topological space (set of points of topological space where well formed formula is true ). If a topological space $(X, \tau)$ and valuation function $\nu :\textsc{Prop}\to \mathcal{P}(X)$ is fixed, there is unique expansion of the valuation function to the set of all well formed formulas in the language $\mathcal{ML}$.

\begin{center}
\begin{tabular}{lcl}
$\nu(\bot )=\emptyset $ & & \\
$x\models p$ & $\mathrm{ iff}$ & $x\in \nu (p)$ \\
$x\models \neg \phi $ & $\mathrm{ iff}$ & $x \not \models \phi$ \\
$x\models \phi \vee \psi$ & $\mathrm{ iff}$ & $x\models \phi $ or $x\models \psi$ \\
$x\models \Diamond \phi $ & $\mathrm{ iff}$ & $\forall O\in \tau, \ x\in O \ \text{implies} \ \exists y\in O \ \text{such that} \ y\models \phi  $\\
\end{tabular}
\end{center}
And hence
\begin{center}
\begin{tabular}{lcl}
$x\models \Box \phi $ & $\mathrm{ iff}$ & $\exists O\in \tau \ \text{ such that } \  x\in O \ \text{ and } \ \forall y\in O\ \  y\models \phi $
\end{tabular}
\end{center}
That is, the regions denoted by the propositional letters are specified in advance by means of a valuation function, and $\vee$, $\neg$,  $\Diamond$ and $\Box$ are interpreted as union, complementation, the closure and  the interior operator respectively.

Each modal formula $\phi$ defines a set of points in a topological model (namely the set of points at which it is true). With a slight overloading of notation, we will sometimes denote this set by $\nu(\phi)$. It is not hard to see that $\nu (\Box \phi)=\mathbb{I} \nu(\phi)$ and $\nu (\Diamond \phi)=\mathbb{C} \nu(\phi)$. A topological space $(X, \tau)$ with fixed valuation $\nu$ is called a \emph{topological model} (\emph{topo-model}) of the modal language $\mathcal{ML}$ and is denoted by $\mathfrak{M}=(X,\tau , \nu)$.

Occasionally we will also equivalently talk about the set $\nu(\phi)$ in terms of its characteristic function, i.~e. view it as a map $\nu(\phi):X\to\{\textbf{true},\textbf{false}\}$ which sends $x\in X$ to the truth value of $x\models\phi$.

For a subset $A\subset X$ if for all $x\in A$ holds $\mathfrak{M},x\models \phi$, then   we shortly write $\mathfrak{M},A\models \phi$. Further, $\mathfrak{M}\models \phi$ ($\phi$ is valid in $\mathfrak{M}$) means that $\mathfrak{M},x\models \phi$ for all $x\in X$. We write $X\models \phi$ ($\phi$ is valid in $X$) when $(X,\nu)\models \phi$ for any valuation $\nu$. If $\mathbf{K}$ is a class of models (topomodels in our case) we write $\mathbf{K}\models \phi$ when $X\models \phi$ for each $X\in \mathbf{K}$. By $\mathsf{Log}(\mathbf{K})$ we denote the set of all modal formulas valid in all members $X\in \mathbf{K}$. In case $\mathbf{K}$ consists of a single member $X$ we may write $\mathsf{Log}(X)$ to denote the modal logic of $X$.

The modal logic $\mathbf{S4}$ of Lewis is defined as follows:\\
The modal logic $\mathbf{S4}$ is the smallest set of modal formulas which contains all the classical tautologies, the following axioms:\\

\begin{tabular}{llll}
 Axiom & Diamond notations &  or & Box notations \\
(N)& $\neg \Diamond \bot $ & & $\Box \top$ \\
(T)& $p\rightarrow \Diamond p$ & & $\Box p \rightarrow p$ \\
(R)& $\Diamond (p \vee q) \leftrightarrow (\Diamond p \vee \Diamond q )$ & & $\Box (p \wedge q) \leftrightarrow (\Box p \wedge \Box q )$\\
(4)& $\Diamond \Diamond p \rightarrow \Diamond p $ & & $\Box p \rightarrow \Box\Box p$
\end{tabular}

\ 

\noindent and is closed under the rules of modus ponens and necessitation (from $\phi$ derive $\Box \phi$ ), and under uniform substitution.\\

It is well known that the logic $\mathbf{S4}$ is characterized by reflexive transitive Kripke frames. As we mentioned such Kripke frames can naturally be seen as Alexandroff topological spaces. Moreover these are the classical results of McKinsey and Tarski which state in particular, that the modal logic $\mathbf{S4}$ is sound and complete w.r.t. the class of all topological models, as well as $\mathsf{Log}(\mathbb{R}^n)=\mathbf{S4}$ for any Euclidean space $\mathbb{R}^n$ \cite{mck:mckinsey44}.

\subsection{Topological bisimulations}
As a main semantical tool the paper employs the following generalization of the notion of \emph{bisimulation} between Kripke models. An interested in more extensive reference reader may consult with textbook like \cite{spatlog:AielloBenthem07}.

  Consider topological models $\mathfrak{M}=(X,\nu)$ and $\mathfrak{M}'=(X',\nu')$. A non-empty relation $Z\subseteq X\times X'$ is a \emph{topo-bisimulation} between $\mathfrak{M}$ and $\mathfrak{M}'$ if the following conditions are met for all $x\in X$ and $x'\in X'$:\\
 \begin{tabbing}
 \textbf{Atom}\= \hspace{0.5cm} if $xZx'$ then $x\in\nu(p)$ iff $x'\in\nu'(p)$ for all $p\in\textsc{Prop}$\kill
  \textbf{Atom}\> \hspace{0.5cm} if $xZx'$ then $x\in\nu(p)$ iff $x'\in\nu'(p)$ for all $p\in\textsc{Prop}$.\\
  \textbf{Zig}\> \hspace{0.5cm} For all $O\in\tau$, $Z[O]\in\tau'$.\\
  \textbf{Zag}\> \hspace{0.5cm} For all $O'\in\tau'$, $Z^{-1}[O']\in\tau$.\\
\end{tabbing}

Here $Z[O]=\{x'|\exists x\in O(xZx')\}$ denotes the image of $O$ under $Z$. The preimage $Z^{-1}[O']$ is defined analogously.

Topo-bisimulations are closely linked with the notion of modal equivalence. Evidence for this comes from the following result \cite{bent:AielloBenthem02}:
\begin{theorem}
\label{t:bisimimplmodeq}
Let $\mathfrak{M}_1=(X_1,\nu_1)$ and  $\mathfrak{M}_2=(X_2,\nu_2)$ be two topo-models, and $x\in X_1$ and $x'\in X_2$ be two topo-bisimilar points. Then for each modal formula $\phi$, condition $\mathfrak{M}_1,x\models \phi$ holds if and only if $\mathfrak{M}_2,x'\models \phi$ holds. That is, modal formulas are invariant under topo-bisimulations.
\end{theorem}

There is the similarity between interior maps and bisimulations. Both require images and pre-images of opens to be open. Indeed this similarity can be utilized to show that onto interior maps between topological spaces preserve modal validity. Here we present a slightly more general result which will be put to use later. Let $f: X\pto Y$ denote a partial map from $X$ to $Y$. When the domain of $f$ is an open subset of $X$, and in addition $f$ is an interior mapping from it's domain to $Y$, we call $f$ an interior partial map.

\begin{proposition}\label{t:partial_interior_map}
Let $X$ and $Y$ be topological spaces and let $f:X\pto Y$ be an onto partial interior map. Then for an arbitrary modal formula $\phi$ we have $Y\models\phi$ whenever $X\models\phi$.
\end{proposition}
\begin{proof}
Let $U\subseteq X$ be open and let $i:U\to X$ be the identity map. We proceed by contraposition. Suppose $Y\not\models\phi$ for some modal formula $\phi$. Then there is a valuation $\nu$ on $Y$ and a point $y\in Y$ such that $Y,\nu,y\models\neg\phi$. Consider a valuation $\nu_U$ on $U$ defined by $\nu_U(p)=f^{-1}(\nu(p))$. Moreover, consider the same valuation on $X$, that is put $\nu_X(p)=f^{-1}(\nu(p))=i(\nu_U(p))$. Now pick any point $x\in X$ such that $f(x)=y$. It is straightforward to check that the pointed model $(Y,\nu,y)$ is topo-bisimlar to $(U,\nu_U,x)$ which in turn is topo-bisimilar to $(X,\nu_X, x)$. The topo-bisimulations are provided by the graphs of the maps $f$ and $i$, respectively. It readily follows that $(Y,\nu,y)$ is topo-bisimilar to $(X,\nu_X, x)$. Hence, by Theorem~\ref{t:bisimimplmodeq} we conclude that $X,\nu_X, x\models\neg\phi$ and thus $X\not\models\phi$.
\end{proof}

It follows that if $Y$ is a partial interior image of $X$, then ${\sf Log}(X)\subseteq {\sf Log}(Y)$.

\subsection{General spaces}

It is well known that neither Kripke semantics, nor topological semantics is fully adequate for modal logic. In case of Kripke semantics, a fully adequate generalization is provided by general frames, where valuations are restricted to modal subalgebras of the powerset algebra of a Kripke frame. A similar approach has been introduced for topological spaces in \cite{spatlog:AielloBenthem07}: they consider \emph{general spaces}, which are topological spaces together with a fixed collection of subsets that is closed under set-theoretic operations as well as under operation of topological closure. In the main part of the paper we are interested in specific valuations in Euclidean spaces that only take as values geometrically simple sets, namely in the Boolean algebra with closure generated by planar polygons.

\begin{definition}
\label{gentopmodel}
A \emph{general topological model} is a triple $(X,A,\nu)$, where 
\begin{itemize}
\item{$X=(X,\tau)$ is a
topological space;}
\item{$A$ is a modal subalgebra of $\mathcal{P}(X)$ (i.~e. subset of $\mathcal{P}(X)$ closed under Boolean operations and under the operation of closure $\mathbb{C}$);}
\item{$\nu:~Prop\to A$ (valuation) is a map which sends propositional letters to specific subsets of $X$, namely the elements of $A$.}
\end{itemize}
We refer to $(X,A)$ as \emph{general space} and to elements of $A$ as \emph{admissible sets}.  \\
  A morphism between general spaces $\mathfrak{M}_1=(X_1,\tau_1, A_1)$ and $\mathfrak{M}_2=(X_2,\tau_2,A_2)$ is an interior map $f:X_1 \to X_2$ such that $f^{-1}(U)\in A_1$ for each $U\in A_2$.
\end{definition}

The notion of validity for general spaces is defined as expected. It is a simple exercise to check that Theorem~\ref{t:partial_interior_map} extends to the case of the partial onto interior maps between general spaces.

\subsection{Ultraproducts and first order definability}

\begin{definition}
A class \textsf{K} of models for fixed first-order language $\mathcal{L}$ is \emph{defined by} a set $\Sigma$ of $\mathcal{L}$-sentences if every model for the language is in \textsf{K} iff it is a
model for $\Sigma$. A class of models is \emph{elementary} if it is defined by some set of first-order sentences.
\end{definition}
To deal with the notion of first order definability we employ the notion of \emph{ultraproduct}. The purpose of that will be shown below. Suppose $U$ is an ultrafilter over a nonempty set $I$, and $(X_i)$ is a family of nonempty sets indexed by $I$. Let $P=\prod_{i\in I} X_i$ be the cartesian product of the family. An element of $P$ are functions from $I$ to the disjoint union of the $X_i$'s, such that $f(i)\in X_i$ for each $i\in I$. Given two elements $f,g\in P$  we say that $f$ and $g$ are $U$-equivalent, (and denote it by $f\sim_{U} g $) if $\{i\in I| f(i)=g(i)\}\in U$. The relation $\sim_{U}$ defined in this way is an equivalence relation.
\begin{definition}
Let $f_U$ denote the equivalence class of $f$ modulo $\sim_U$. The \emph{ultraproduct of the sets $X_i$ modulo $U$} is the set of all equivalence classes of $\sim_U$. It is denoted by $\prod_{U}X_i$
\end{definition}
One can easily expand the above definition of ultraproduct to models:

\begin{definition}
Fix a first-order language $\mathcal{L}$, and let $\mathfrak{M}_i$ be models of $\mathcal{L}$ indexed by the set $I$. The $\mathit{ultraproduct}$  $\prod_{U}\mathfrak{M}_i$ modulo $U$ is the model described as follows:

\begin{enumerate}[(i)]
\item{The underlying set of $\prod_{U}\mathfrak{M}_i$ is a set $\prod_{U}X_i$ }, where $X_i$ is the underlying set of $\mathfrak{M}_i$.
\item{Let $R$ be an $n$-place relation symbol, and $R_i$ its interepretation in the model $\mathfrak{M}_i$ The relation $R_U$ in $\prod_{U}\mathfrak{M}_i$ is given by
$$R_{U}(f_{U}^{1},\dots , f_{U}^{n}) ~\text{~iff~} ~\{i\in I| ~R_{i}\big(f^{1}(i),\dots ,f^{n}(i) \big) \}\in U.$$}
\item{Let $F$ be a $m$-place function symbol, and $F_i$ its interepretation in the model $\mathfrak{M}_i$ The function $F_U$ in $\prod_{U}\mathfrak{M}_i$ is given by
$$F_{U}(f_{U}^{1},\dots , f_{U}^{n})= \{\big( i,F_{i}(f^{1}(i),\dots ,f^{n}(i))\big)|~i\in I\}_U.$$}
\item{Let $c$ be a constant, and $a_i$ its interpretation in $\mathfrak{M}_i$. Then $c$ is interpreted by the element $c'\in \prod_{U}X_{i}$ where $c'=\{(i,a_i)|~i\in I\}$ }
\end{enumerate}
In the case where all the models are the same, say $\mathfrak{M}_{i}=\mathfrak{M}$, we speak of \emph{ultrapower} of $\mathfrak{M}$ modulo $U$.
 \end{definition}
 The following proposition connects the notions of ultraproduct and first-order definability  
\begin{theorem}{(\cite{changkeisler73} Theorem 6.1.16)}{\label{fod}}
A class of models \textsf{K} is defined by means of a set of first-order sentences if and only if it is closed under isomorphisms and ultraproducts, while its complement
is closed under ultrapowers.
\end{theorem}
After this preliminary part we are ready to deal with  the main matter of the paper.
\section{The Euclidean polytopal spaces, the Logic $\mathbf{PL_2}$ and the finite model property}\label{fmp}
We are interested in specific general spaces defined over Euclidean spaces $\mathbb{R}^n$. Let $P_n$ be the boolean algebra of the $n$-dimensional polytopes in $\mathbb{R}^n$. By a \emph{polytope} we mean a finite union of subsets in $\mathbb{R}^n$ which are solutions of a system of linear inequalities.
Alternatively, simple polytopes can be described as sets that are intersections of finitely many hyperplanes in $\mathbb{R}^n$, and the polytopes are the finite unions of those. To be more precise, a \emph{polytope} is any subset of $\mathbb{R}^n$ of the form $P=\{\bar{x}\mid \bigvee\bigwedge (\ell_i(\bar{x})\bowtie a_i)\}$ where $\ell_i$ are linear forms on $\mathbb{R}^n$, where $a_i$ are real numbers, $\bowtie$ denotes any of the inequality symbols $\ge,>, \le, <$, while $\bigvee$ and $\bigwedge$ denote finite disjunction and finite conjunction. The sets of the form $P=\{\bar{x}\mid \bigwedge (\ell_i(\bar{x})\bowtie a_i)\}$ we call \emph{simple polytopes}.

Then  it is clear from the definition that the set $P_n$ of all the $n$-dimensional polytopes forms a boolean subalgebra of the powerset of $\mathbb{R}^n$ (note that the negation of an inequality is again an inequality). Moreover we have

\begin{proposition}\label{t:cloclo}
The boolean algebra $P_n$ is a modal subalgebra of the powerset of $\mathbb{R}^n$ equipped with the closure/interior operators for the Euclidean topology on $\mathbb{R}^n$.
\end{proposition}
\begin{proof}
To show that $P_n$ is closed under the closure operator, first note that the closure operator distributes over finite unions. Since every polygon is a finite union of simple polygons, it suffices to point out that given a nonempty simple polygon defined by a finite conjunction of inequalities, its closure is just the closed simple polygon obtained by turning all strict inequalities into non-strict ones.
\end{proof}

\begin{definition}
Let $\mathfrak P_n=(\mathbb{R}^n, P_n)$ be the general space defined by polytopes. We call such a space the \emph{$n$-dimensional Euclidean polytopal space}.  The modal logic $\mathbf{PL_n}$ of the $n$-dimensional Euclidean polytopal  space is defined to be the set of all modal formulas which are valid on $\mathfrak P_n$.
$$\mathbf{PL_n}\coloneqq {\sf Log}(\mathfrak P_n)$$
\end{definition}

In this paper we concentrate on 2-dimensional polytopal modal logic $\mathbf{PL_2}$. We call this logic the polygonal modal logic for brevity and the corresponding general space $\mathfrak P_2=(\mathbb{R}^2, P_2)$ we call \emph{polygonal plane}. 

The admissible sets in $\mathfrak P_2$  are finite unions of \emph{generalized planar polygons}, where under a generalized planar polygon we understand a (possibly unbounded) region in the plane which is an intersection of finitely many (closed or open) half-planes. It is clear that any point, line, ray or segment also falls under this definition, as do triangles, pentagons and $n$-gons in general. 

However for our modal-logical purposes, we may restrict attention to usual bounded polygons, line segments, points, their complements and finite unions of these. The reason of this is as follows: First of all note, that these bounded and co-bounded polygons together form a boolean algebra themselves. Now suppose a modal formula $\phi$ is satisfiable at a point $x\in\mathbb{R}^2$ with an admissible valuation $\nu :\textsc{Prop} \to P_n$. If we take any open triangle $T$ containing $x$ and change the valuation so that $\nu'(p)=T\cap\nu(p)$. Then $\mathbb{R}^2,\nu',x\models \phi$ and $\nu'$ is now admitting only bounded polygons. On other hand is clear that $(\mathbb{R}^2,\nu,x)$ and $(\mathbb{R}^2,\nu',x)$ are topo-bisimilar to each other and by Proposition \ref{t:bisimimplmodeq} we may restrict our attention to bounded polygons.

Now we are going to investigate what kind a finite Kripke frame $(X,R)$ should be to satisfy condition $(X,R)\models\mathbf{PL_2}$. It is easy to understand, that if we are able to investigate finite Kripke frames of $\mathbf{PL_2}$, then it facilitates the task of study of $\mathbf{PL_2}$ itself. Hence to investigate the question of finite Kripke models of $\mathbf{PL_2}$, we need to look at finite interior (open and continuous) images of the polygonal plane. Here we should note that we need to preserve formulas and hence we need a mapping from \emph{polygonal} plane induce a homomorphism of modal algebras.  Hence we have to investigate finite interior images of the polygonal plane such that the pre-image of each point is a finite non-empty union of planar polygons. 

We consider the class of Kripke frames we call \emph{crown frames} (not to be confused with the \emph{crown graph}).
\begin{definition}\label{saw}
A \emph{crown  frame} $\SAW_n$ is a frame $(S_n,Q_n)$ such that  $S_n=\{ r, s_{1},\cdots , s_{2n} \} $ and $Q_n$ is defined as follows:	
\begin{center}
\begin{tabular}{ll}
$rQ_n r $ &\hspace{0.5cm}where $r$ is a root,\\
$s_i Q_n  s_i$ & \hspace{0.5cm}for all  $s_i\in S_n$,\\
$rQ_n s_i$ & \hspace{0.5cm}for all  $s_i\in S_n$,\\
$s_i Q_n s_j$ & \hspace{0.5cm}when $i<2n$ is even and $j=i-1,i+1$ ,\\
$s_{2n}Q_n s_1 $ & \hspace{0.5cm}and  $s_{2n}Q_n s_{2n-1} $.  \\
\end{tabular}
\end{center}
\end{definition}

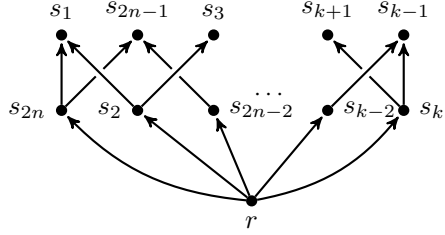
\begin{figure}[H]

\begin{center}
\begin{tikzpicture} [->, >=stealth', shorten >=1pt,auto, node distance=1cm,
 thick,
 upa node/.style={circle,inner sep=1.2pt, fill=black, draw, font=\sffamily \Large \bfseries, label={$s_1$}},
 upb node/.style={circle,inner sep=1.2pt, fill=black, draw, font=\sffamily \Large \bfseries, label={$s_{2n-1}$}},
 upc node/.style={circle,inner sep=1.2pt, fill=black, draw, font=\sffamily \Large \bfseries, label={$s_3$}},
 upd node/.style={circle,inner sep=1.2pt, fill=black, draw, font=\sffamily \Large \bfseries, label={$s_{k+1}$}},
 upe node/.style={circle,inner sep=1.2pt, fill=black, draw, font=\sffamily \Large \bfseries, label={$s_{k-1}$}},
 da node/.style={circle,inner sep=1.2pt, fill=black, draw, font=\sffamily \Large \bfseries, label={left:$s_{2n}$}},
 db node/.style={circle,inner sep=1.2pt, fill=black, draw, font=\sffamily \Large \bfseries, label={left:$s_2$}},
 dc node/.style={circle,inner sep=1.2pt, fill=black, draw, font=\sffamily \Large \bfseries, label={right:$s_{2n-2}$}},
 dd node/.style={circle,inner sep=1.2pt, fill=black, draw, font=\sffamily \Large \bfseries, label={right:$s_{k-2}$}},
 de node/.style={circle,inner sep=1.2pt, fill=black, draw, font=\sffamily \Large \bfseries, label={right:$s_{k}$}},
 r node/.style={circle,inner sep=1.2pt, fill=black, draw, font=\sffamily \Large \bfseries, label={below:$r$}}]

  \node[upa node] (1) {};
  \node[upb node] (2) [right of=1] {};
  \node[upc node] (3) [right of=2] {};
  \node[upd node] (4) [right of=3,xshift=0.5cm] {};
  \node[upe node] (5) [right of=4] {};
  \node[da node] (6) [below of=1] {};
  \node[db node] (7) [below of=2] {};
  \node[dc node] (8) [below of=3] {};
  \node[dd node] (9) [below of=4] {};
  \node[de node] (10) [below of=5] {};
  \node[r node] (r) [below of=8,xshift=0.5cm, yshift=-0.2cm] {};

  \path[every node/.style={font=\sffamily\small}]

    (6) edge [] node[] {} (1)
        edge [] node[] {} (2)
    (8) edge [] node[] {} (2)
        edge [draw=none] node[] {$\cdots$} (9)
    (7) edge [-,draw=white, line width=4pt] node[] {} (1)
        edge [] node[] {} (1)
       edge [-,draw=white, line width=4pt] node[] {} (3)
        edge [] node[] {} (3)
    (10)edge [] node[] {} (4)
        edge [] node[] {} (5)
    (9) edge [-,draw=white, line width=4pt] node[] {} (5)
         edge [] node[] {} (5)
    (r) edge [] node[] {} (7)
        edge [] node[] {} (8)
        edge [] node[] {} (9);

       \draw[->] (r) to [bend right=20] node[] {} (10);
       \draw[->] (r) to [bend left=20] node[] {} (6);

\end{tikzpicture}
\par\end{center}

  \centering
 \caption{Crown frame}
 \label{sawfig}
\end{figure}

A general crown frame is depicted on figure~\ref{sawfig}.

The following proposition states that crown frames are typical examples of finite interior images of the polygonal plane $\mathfrak P_n$.

\begin{proposition}\label{t:fmp easy direction}
For any number $n\in \mathbb{N}$ the crown frame $\SAW_n$ is an interior image of the polygonal plane $\mathfrak{P}_2$.\footnote{By virtue of mentioned above correspondence between quasi-ordered sets, i.e. reflexive and transitive Kripke frames, here we consider a crown frame as an Aleksandroff topological space obtained by declaring its upwards closed subsets open.}
\end{proposition}
\begin{proof}

Suppose a number $n$ for $\SAW_n$ is fixed. Consider any point $x$ and arbitrary distinct rays $l_1,\dots, l_n$ emanating from $x$ and enumerated in the clockwise direction.
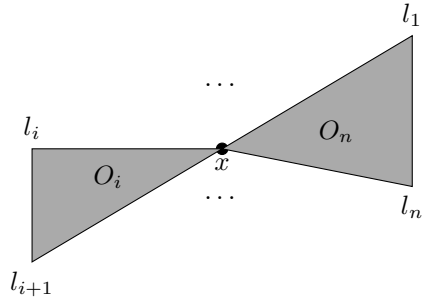
\begin{center}
\begin{tikzpicture}
\draw[black,fill=black] (0,0) circle (.5ex);
\coordinate[label={[below]$x$}] (Origin)   at (0,0) circle (.5ex);
\coordinate [label={[above]$l_1$}](Lone) at (2.5,1.5);
\coordinate [label={[below]$l_n$}] (Ln) at (2.5,-0.5);
\coordinate [label={[above]$l_i$}](Li) at (-2.5,0);
\coordinate [label={[below]$l_{i+1}$}](Liplusone) at (-2.5,-1.5);
\coordinate (corneronen) at (2.5,1.5);
\coordinate (corneriiplusone) at (-2.5,-1.5);
\coordinate [label={[below]$\dots$}](Dotsup) at (0,1);
\coordinate [label={[below]$\dots$}](Dotsdown) at (0,-0.5);

\draw [-,fill={rgb:black,1;white,2}] (Origin) -- (Lone) -- (corneronen) -- (Ln) -- cycle;
\draw [-,fill={rgb:black,1;white,2}] (Origin) -- (Li) -- (corneriiplusone) -- (Liplusone) -- cycle;

\coordinate [label={[below]$O_n$}](On) at (1.5,0.5);
\coordinate [label={[below]$O_i$}](Oi) at (-1.5,-0.1);

\end{tikzpicture}
\end{center}
Let $O_i$ denote the open region between $l_i$ and $l_{i+1}$, where $i\in\{1,\dots,n-1\}$ and let $O_n$ denote the open region between $l_n$ and $l_1$. \\
Define the map $f:\mathbb{R}^2\to S_n$ by putting $f(x)=r$, $f(l_i)=s_{2i}$ and $f(O_i)=s_{2i-1}$. It is straightforward to check that $f$ is an interior map. We only expand one of the less obvious steps in the proof. Suppose $U\subseteq \mathbb{R}^2$ is open and suppose $s_{2k}\in f(U)$. Then there exists a point $y\in l_k$ such that $y\in U$ and $f(y)=s_{2k}$. Since $U$ is open, if it contains a point from the ray $l_k$, it will intersect both regions neighbouring $l_k$, namely $O_{k-1}$ and $O_{k+1}$. It follows that $f(U)$ will contain both of the points from $O\subseteq S_n$ that $s_{2k}$ is related to. This shows that the image of $U$ under $f$ is an upset. Hence $f$ is an open map.
\end{proof}

The intuition employed in the proof of the above proposition is as follows: if a point $x$ is `close' to a set $A$ in the plane, then the image of $x$ is related to the image of $A$ under an interior map. To make it clear, this is encoded in the following equivalence: $x\in\mathbb{C} f^{-1}(u)\ \ \mathrm{iff}\ \ f(x) R u$, which is similar to the p-morphism condition for Kripke frames. In other words, interior maps preserve the `closeness' relation between points and sets: $x\in\mathbb{C}A\ \ \mathrm{iff}\ \ f(x)\in R^{-1}[f(A)]$.

In light of the Theorem~\ref{t:partial_interior_map}, any modal formula satisfiable on one of the crown frames is also satisfiable on the polygonal plane. Our next aim is to show the converse.
Actually we are going to show that the logic of the class of all crown frames coincides with $\mathbf{PL}_2$.
\begin{theorem}
\label{t:fmp hard direction}
Let $\phi$ be satisfiable on a polygonal plane. Then $\phi$ is satisfiable on one of the crown frames.
\end{theorem}
\begin{proof}

Let $\nu$ be a valuation and $x$ be a point such that $\mathfrak P_2,\nu,x\models\phi$.

Our strategy will be to find a small enough open neighborhood $U$ around $x$ such that the partial interior map could be built from $U$ onto one of the crown frames in such a way that the pre-images of possible worlds 
from the crown frame have constant valuation for propositional variables occurring in $\phi$.

Suppose $\phi$ depends on propositional variables $p_1,\dots,p_k$. It is clear that the truth of $\phi$ will not be affected if we assume that all the other propositional letters are mapped to empty set. Let $A_i\Coloneqq\nu(p_i)$ for $i\in\{1,\dots,k\}$. Then each $A_i$ is a finite union of simple polygons. Let $S$ be the collection of all the simple polygons occurring in the $A_i$s. Let $E$ be the collection of all lines occurring as an edge of one of the simple polygons in $S$ (in case there is only one such $A_i$, then draw any segment with endpoint $x$). It is obvious that $E$ is finite. Furthermore, we observe that for any segment $I$ on the plane, if the endpoints of the segment differ on the valuation of a propositional letter $p_i$, then the segment $I$ must intersect with one from $E$, namely the one that is represented as a border of $A_i$ which $I$ must cross in order to change valuation from one endpoint to the other.

Now, for each line in $E$, calculate the distance from $x$ to that line and to its endpoints (if it have such). This will produce a finite number of non-negative real numbers. Let $\alpha$ be the least \emph{positive} number thus obtained and let $B=B(x,{\frac{\alpha}{2} })$ be the open ball with the center at $x$ and the radius $\frac{\alpha}{2} $.
It is straightforward that only the lines from $E$ that pass through $x$ (or have it as an endpoint) will intersect with $B$. Let us label the intersection points of lines from $E$ with the border of $B$ in a clockwise direction as $x_1, x_2,\dots,x_m$ with $m\leq k$. Let $l_i$ denote the open segments $(x,x_i)$ and let $O_i$ denote the open sectors of $B$ bounded by $l_i$ and $l_{i+1}$ for $i\in\{1, \dots, m-1\}$. Let $O_m$ be the remaining open sector defined by $l_m$ and $l_1$. Then $B$ breaks down into the sets $\{x\}$, $l_i$ and $O_i$.
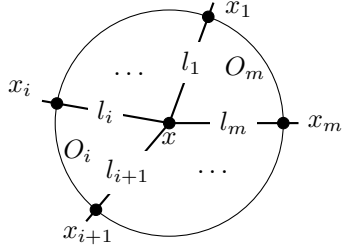
\begin{center}
\begin{tikzpicture}
\draw[black,fill=black] (0,0) circle (.5ex);
\coordinate[label={[below]$x$}] (Origin)   at (0,0) circle (.5ex);
\draw (0,0) circle (1.5cm);
\draw [-,thick]  (0,0) -- node[fill=white, align = center]{$l_m$} (xyz polar cs:angle=0,radius=1.7);
\draw [-,thick]  (0,0) -- node[fill=white, align = center]{$l_1$} (xyz polar cs:angle=70,radius=1.7);
\coordinate[label={[below]$\dots$}] (Dots) at (-0.5,0.8);
\draw [-,thick]  (0,0) -- node[fill=white, align = center]{$l_i$} (xyz polar cs:angle=170,radius=1.7);
\draw [-,thick]  (0,0) -- node[fill=white, align = center]{$l_{i+1}$} (xyz polar cs:angle=230,radius=1.7);
\coordinate[label={[below]$\dots$}] (Dots) at (0.6,-0.5);
\coordinate[label={[below]$O_m$}] (Dots) at (1,1);
\coordinate[label={[below]$O_i$}] (Dots) at (-1.2,-0.1);

\coordinate[label={[right]$x_1$}] (x_1)   at (0.61303021,1.50953893) circle (.5ex);
\draw[black,fill=black] (0.51303021,1.40953893) circle (.5ex);

\coordinate[label={[right]$x_m$}] (x_m)   at (1.7,0) circle (.5ex);
\draw[black,fill=black] (1.5,0) circle (.5ex);

\coordinate[label={[left]$x_i$}] (x_i)   at (-1.67721162951,0.46047226649) circle (.5ex);
\draw[black,fill=black] (-1.47721162951,0.26047226649) circle (.5ex);

\coordinate[label={[below]$x_{i+1}$}] (x_i+1)   at (-1.06418141452,-1.24906666466) circle (.5ex);
\draw[black,fill=black] (-0.96418141452,-1.14906666466) circle (.5ex);

\end{tikzpicture}
\end{center}
We define the map $f$ from $B$ onto $\SAW_m$ in the obvious way (see the proof of Proposition~\ref{t:fmp easy direction} from above). This map is easily seen to be an interior map. We claim more, namely that the valuation of the propositional variables $p_j$ occurring in $\phi$ is constant on each of the sets $E_i$ and $O_i$. Indeed, take any two points $y,z\in O_i$ for some $i\in\{1,\dots,m\}$ and suppose $y\in A_j\not\ni z$ for some $j\in\{1,\dots,k\}$. But then the whole closed segment $[y,z]$ falls inside $O_i$ since the latter is convex. On the other hand, the endpoints of this segment differ on the valuation of $p_j$, which means that the segment $[y,z]$ must cross a member of $E$, which is impossible. A similar argument shows that no two points can disagree on the valuation of $p_j$ along any of the $E_i$.

It follows that the valuation $\mu$ defined by putting $\mu(p)=f(\nu(p))$ for each $p$ is such that

$$f(x)\in\mu(p)\ \ \textrm{iff}\ \ x\in\nu(p).$$

Consequently, the graph of $f$ is a topo-bisimulation between $(\mathbb{R}^2,\nu,x)$ and $(\SAW_m,\mu, r)$. It follows that $\phi$ is satisfiable on $\SAW_m$.

\end{proof}
Hence we have proved completeness of the logic with respect to class of finite crown frames. So the next corollary immediately follows:
\begin{corollary}
The logic $\mathbf{PL}_2$ is determined by the class of finite crown frames. Hence this logic has fmp.
\end{corollary}

\begin{proof}
Just putting together Proposition~\ref{t:fmp easy direction} and Theorem~\ref{t:fmp hard direction}.
\end{proof}

Therefore, the task of subsequent investigation of $\mathbf{PL}_2$ can be reduced to the setting of Kripke semantics.

\section{Axiomatization}\label{axiomatization}
As we proved fmp for $\mathbf{PL}_2$, the logic under consideration, our next task is to axiomatize  the logic. To do that we use Jankov-Fine formulas for finite rooted frames that are \emph{not} $\mathbf{PL}_2$-frames. First of all, we describe the five simplest frames that falsify $\mathbf{PL}_2$. These are depicted below.

\begin{figure}[H]
\label{forbidden}
\begin{center}
\begin{tikzpicture}[->,>=stealth',shorten >=1pt,auto,node distance=1cm,
  thick,main node/.style={circle,inner sep=1.2pt,fill=black,draw,font=\sffamily\Large\bfseries},every fit/.style={ellipse,draw,inner sep=1.2pt}]

  \node[main node] (1) [] {};
  \node[main node] (2) [right of=1] {};

  \node[draw,fit= (1) (2),minimum height=1cm] {} ;

  \node[] (caption) [below right of=1,xshift=-0.2cm,yshift=-0.2cm,draw=none,fill=none] {$\mathfrak B_1$};

\end{tikzpicture}
\hspace{0.4cm}
\begin{tikzpicture}[->,>=stealth',shorten >=1pt,auto,node distance=2cm,
  thick,main node/.style={circle,inner sep=1.2pt,fill=black,draw,font=\sffamily\Large\bfseries},every fit/.style={ellipse,draw,inner sep=1.2pt}]

  \node[main node] (1)[]{};
  \node[main node] (2) [below left of=1,xshift=0.95cm] {};
  \node[] (3) [below of=1,yshift=0.95cm,draw=none,fill=none] {};
  \node[main node] (4) [right of=2,xshift=-1cm] {};

  \path[every node/.style={font=\sffamily\small}]

    (3) edge node [] {} (1);

  \node[draw,fit= (2) (4),minimum height=1cm] {} ;

  \node[] (caption) [below of=3,yshift=0.7cm,draw=none,fill=none] {$\mathfrak B_2$};

\end{tikzpicture}
\hspace{0.1cm}
\begin{tikzpicture}[->,>=stealth',shorten >=1pt,auto,node distance=1cm,
  thick,main node/.style={circle,inner sep=1.2pt,fill=black,draw,font=\sffamily\Large\bfseries}]

  \node[main node] (1)[] {};
  \node[main node] (2) [right of=1] {};
  \node[main node] (3) [right of=2] {};
  \node[main node] (4) [below of=2] {};

  \path[every node/.style={font=\sffamily\small}]

    (4) edge [] node[] {} (1)
         edge [] node[] {} (2)
         edge [] node[] {} (3);

  \node[] (caption) [below of=4,yshift=0.4cm,draw=none,fill=none] {$\mathfrak B_3$};S

\end{tikzpicture}
\hspace{0.3cm}
\begin{tikzpicture}[->,>=stealth',shorten >=1pt,auto,node distance=1cm,
  thick,main node/.style={circle,inner sep=1.2pt,fill=black,draw,font=\sffamily\Large\bfseries}]

  \node[main node] (1) []{};
  \node[main node] (2) [below of=1] {};
  \node[main node] (3) [below of=2] {};
  \node[main node] (4) [below of=3] {};

  \path[every node/.style={font=\sffamily\small}]

    (2) edge [] node[] {} (1)
    (3) edge [] node[] {} (2)
    (4) edge [] node[] {} (3);

  \node[] (caption) [below of=4,yshift=0.4cm,draw=none,fill=none] {$\mathfrak B_4$};

\end{tikzpicture}
\hspace{0.4cm}
\begin{tikzpicture}[->,>=stealth',shorten >=1pt,auto,node distance=1cm,
  thick,main node/.style={circle,inner sep=1.2pt,fill=black,draw,font=\sffamily\Large\bfseries}]

  \node[main node] (1)[] {};
  \node[main node] (2) [below of=1] {};
  \node[main node] (3) [below right of=2,xshift=-0.2cm] {};
  \node[main node] (4) [right of=2] {};

  \path[every node/.style={font=\sffamily\small}]

    (2) edge [] node[] {} (1)
    (3) edge [] node[] {} (2)
         edge [] node[] {} (4);

  \node[] (caption) [below of=3,yshift=0.4cm,draw=none,fill=none] {$\mathfrak B_5$};

\end{tikzpicture}
\par\end{center}

  \centering
 \caption{Frames that falsify $\mathbf{PL}_2$}
\end{figure}
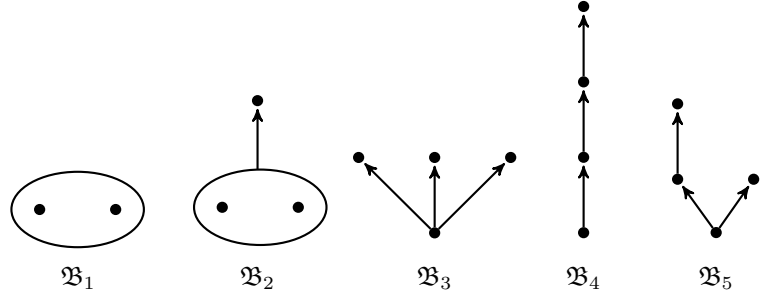

Indeed that these frames cannot support models of $\mathbf{PL}_2$ follows from the following considerations:
\begin{itemize}
\item {- If there would exist a model of $\mathbf{PL}_2$ on $\mathfrak{B}_1$ then there would exist two disjoint finite Boolean combinations of polygons $A_1,A_2$ such that $\mathbb{C}A_{1}= \mathbb{C}A_{2}$, but $A_{1}\neq A_{2}$.}
\item {- If $\mathfrak{B}_2$ would carry a model of $\mathbf{PL}_2$ then there would exist a disjoint finite Boolean combination of polygons $A$ and two disjoint finite Boolean combinations of polygons $A_1,A_2$ such that $\mathbb{C}A_{1}= \mathbb{C}A_{2}$, $A_{1}\neq A_{2}$ and $A_{1}\cup A_{2}\subseteq \mathbb{C}A$.}
\item {- If $\mathfrak{B}_3$ is a model of $\mathbf{PL}_2$ then there would exist a non-empty finite Boolean combination of polygons $A$ which is a subset of closures of three disjoint finite Boolean combinations of polygons $A_1$, $A_2$, $A_3$.}
\item {- If $\mathfrak{B}_4$ is a model of $\mathbf{PL}_2$ then there would exist four disjoint finite non-empty Boolean combinations of polygons $A_1$, $A_2$, $A_3$, $A_4$, such that $$\mathbb{C}A_1\subseteq\mathbb{C}A_2\subseteq\mathbb{C}A_3\subseteq\mathbb{C}A_4.$$ }
\item {- If $\mathfrak{B}_5$ is a model of $\mathbf{PL}_2$ then in the plane there would be two nonempty disjoint open polygons $A$ and $B$ and a polygon $A'\subseteq A$ such that $A'$ is contained in the closure of $A$ and the complement of $A\cup B$ is contained in the closure of $A'$. This would imply that the complement of $A\cup B$ is finite.}
\end{itemize}

All five cases are impossible. We claim that the logic axiomatized by the Jankov-Fine axioms of these five frames coincides with $\mathbf{PL}_2$.
\begin{definition}
Given two Kripke frames $\mathfrak F$ and $\mathfrak G$ we say that $\mathfrak F$ is \emph{subreducible} to $\mathfrak G$ if there exists a generated subframe of $\mathfrak F$ which maps p-morphically onto $\mathfrak G$.
\end{definition}
Note that for $\mathbf{S4}$-frames viewed as topological spaces subreduction is nothing else but the existence of an onto partial interior map. Hence it is not surprising that any formula satisfiable on a frame is also satisfiable on any frame subreducible to it. In fact, for a finite rooted frame, there is a characteristic formula that is satisfiable on precisely the frames subreducible to this frame. Such a formula is called a Jankov-Fine formula of the frame. The precise definition can be found, e.g. in \cite[pp.143]{blackuarn:patrickmaartenyde01}. Let us denote by $\xi(\mathfrak G)$ the Jankov-Fine formula of the frame $\mathfrak G$. Then the following is true:

\begin{proposition}\label{t:Jankov}
If $\xi(\mathfrak G)$ is satisfiable on a frame $\mathfrak F$ then $\mathfrak F$ is subreducible to $\mathfrak G$.
\end{proposition}

The proof can be found in \cite{blackuarn:patrickmaartenyde01}.

Now it follows immediately from the above proposition that the formula $\neg \xi(\mathfrak G)$ is \emph{valid} precisely on those Kripke frames which are \emph{not} subreducible to $\mathfrak G$. The formula $\neg \xi(\mathfrak G)$ is called the \emph{Jankov-Fine axiom} for $\mathfrak G$.

Consider the formula  $\xi=\neg \xi(\mathfrak B_1)\wedge \neg \xi(\mathfrak B_2) \wedge \neg \xi(\mathfrak B_3) \wedge \neg \xi(\mathfrak B_4) \wedge \neg \xi(\mathfrak B_5)$. This formula is valid on a Kripke frame if and only if the frame is not subreducible to any of the frames $\mathfrak B_i$ above. Let us denote by $\Xi$ the logic axiomatized by $\xi$ over $\mathbf{S4}$, i.e. $\Xi\coloneq \mathbf{S4}+\xi$ . We claim that $\Xi=\mathbf{PL}_2$. To show this we first prove that $\Xi$ has the finite model property and then prove that the finite rooted frames for the two logics coincide.

\begin{theorem}\label{t:Xi fmp}
The logic $\Xi$ has the finite model property.
\end{theorem}
\begin{proof}
Note that since $\mathfrak B_4$ is not admitted by $\Xi$, the logic is of finite depth. By Segerberg's Theorem (see, e.g., \cite[Theorem~8.85]{CZ:modlog97}) any logic of finite depth is characterized by its finite frames. It follows that $\Xi$ has the finite model property.
\end{proof}

Since both $\Xi$ and $\mathbf{PL}_2$ are characterized by their finite rooted frames, to show the equality the following two lemmas suffice.

\begin{lemma}\label{l:crowns validate}
Each crown frame validates the axiom $\xi$.
\end{lemma}
\begin{proof}

Let $\SAW_n$ be one of the crown frames.
If (a generated subframe of) $\SAW_n$ is $p$-morphically mapped to a frame containing a cluster then, since $\SAW_n$ is finite and transitive, $\SAW_n$ should already contain a cluster. But it does not. Thus $\SAW_n$ is not subreducible to either one of $\mathfrak B_1$ and $\mathfrak B_2$.
Now assume to the contrary that (a generated subframe of) $\SAW_n$ is
p-morphically mapped by $g$ to the frame $\mathfrak B_3$, the
trident. It is clear that the whole of $\SAW_n$ must be taken, since all the
other generated subframes have no more than three points, while the trident has
four.
Thus $g:\SAW_n\to\mathfrak B_3$ is a p-morphism. Clearly the root is mapped to the root. No other point can be mapped to the root of the trident again because all other points have only two strict successors, while the root of the trident has three. This implies that $\SAW_n$ without the root can be broken down into three disjoint up-sets, which is not the case---it cannot even be broken down into two disjoint up-sets, since it is connected. The contradiction shows that $\SAW_n$ does not subreduce to the trident.
It is easily seen that $\SAW_n$ is not subreducible to $\mathfrak B_4$ because the depth of the p-morphic image of a frame is always less or equal to the depth of the frame.
Finally, to show that $\SAW_n$ is not subreducible to $\mathfrak B_5$ we reason almost exactly as in the case of $\mathfrak B_3$. Suppose, for the sake of contradiction, that $g$ is a subreduction of $\SAW_n$ to $\mathfrak B_5$. Note again, that in that case $g$ must be a reduction, since no point other than the root has three successors in $\SAW_n$. By a similar argument, no point other than the root of $\SAW_n$ can map to the root of  $\mathfrak B_5$. Hence the `upper part' of $\SAW_n$ (everything but the root) breaks down into two disjoint up-sets, which is a contradiction.
Hence, $\SAW_n$ is not subreducible to any of the frames $\mathfrak B_i$, $i=1,\dots,5$. It follows, that $\SAW_n\models\neg\xi(\mathfrak B_i)$ for each $i$. Therefore, $\SAW_n\models\xi$.
\end{proof}

\begin{lemma}\label{l:crowns reduce}
Each rooted finite frame $\mathfrak G$ with  $\mathfrak G\models\xi$ is a subreduction of some crown frame.
\end{lemma}
\begin{proof}
Let $\mathfrak G=(W,R)$ be any rooted frame with $\mathfrak G\models\xi$ and with the root $r$. Then $\mathfrak G$ is not subreducible to
any of the five forbidden frames $\mathfrak B_i$, $i=1,\dots,5$.  We will show that $\mathfrak G$ is either a generated subframe or a p-morphic image of
some crown frame $\SAW_n$.

Let us define $x\overline{R}y$ as $xRy \wedge x\neq y$.
We define a partition the set of states $W$ as follows:
\begin{enumerate}
\item $G_0=\{x\in W\mid xRr\}$
\item $G_1=\{x\in W\mid r\overline{R}x\wedge\nexists y,r\overline{R}y\overline{R}x\}$
\item $G_2=\{x\in W \mid \exists y\in G_1, y\overline{R}x\}$
\end{enumerate}
The sets are disjoint by definition and cover $W$ since the prohibition of the frame $\mathfrak B_4$ ensures that the depth of $\mathfrak G$ is at most $2$.

That $G_0=\{r\}$ follows from the fact that $\mathfrak G$ is not subreducible to
$\mathfrak B_1$ or $\mathfrak B_2$.
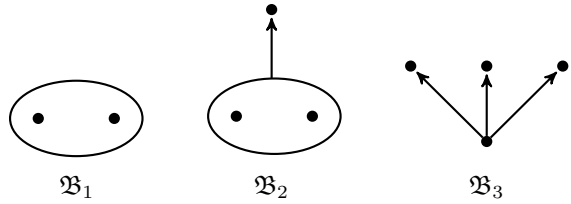
\begin{center}
\begin{tikzpicture}[->,>=stealth',shorten >=1pt,auto,node distance=1cm,
  thick,main node/.style={circle,inner sep=1.2pt,fill=black,draw,font=\sffamily\Large\bfseries},every fit/.style={ellipse,draw,inner sep=1.2pt}]

  \node[main node] (1) [] {};
  \node[main node] (2) [right of=1] {};

  \node[draw,fit= (1) (2),minimum height=1cm] {} ;

  \node[] (caption) [below right of=1,xshift=-0.2cm,yshift=-0.2cm,draw=none,fill=none] {$\mathfrak B_1$};

\end{tikzpicture}
\hspace{0.6cm}
\begin{tikzpicture}[->,>=stealth',shorten >=1pt,auto,node distance=2cm,
  thick,main node/.style={circle,inner sep=1.2pt,fill=black,draw,font=\sffamily\Large\bfseries},every fit/.style={ellipse,draw,inner sep=1.2pt}]

  \node[main node] (1)[]{};
  \node[main node] (2) [below left of=1,xshift=0.95cm] {};
  \node[] (3) [below of=1,yshift=0.95cm,draw=none,fill=none] {};
  \node[main node] (4) [right of=2,xshift=-1cm] {};

  \path[every node/.style={font=\sffamily\small}]

    (3) edge node [] {} (1);

  \node[draw,fit= (2) (4),minimum height=1cm] {} ;

  \node[] (caption) [below of=3,yshift=0.7cm,draw=none,fill=none] {$\mathfrak B_2$};

\end{tikzpicture}
\hspace{0.6cm}
\begin{tikzpicture}[->,>=stealth',shorten >=1pt,auto,node distance=1cm,
  thick,main node/.style={circle,inner sep=1.2pt,fill=black,draw,font=\sffamily\Large\bfseries}]

  \node[main node] (1)[] {};
  \node[main node] (2) [right of=1] {};
  \node[main node] (3) [right of=2] {};
  \node[main node] (4) [below of=2] {};

  \path[every node/.style={font=\sffamily\small}]

    (4) edge [] node[] {} (1)
         edge [] node[] {} (2)
         edge [] node[] {} (3);

  \node[] (caption) [below of=4,yshift=0.4cm,draw=none,fill=none] {$\mathfrak B_3$};S

\end{tikzpicture}
\end{center}
In case $G_2=\emptyset$, since $\mathfrak G$ is not reducible to the trident frame $\mathfrak B_3$, we get $|G_1|\leq 2$ and it is clear that $\mathfrak G$ is (isomorphic to) a generated subframe of either $\SAW_1$ or $\SAW_2$.

\begin{center}
\begin{tikzpicture}[->,>=stealth',shorten >=1pt,auto,node distance=1cm,
  thick,main node/.style={circle,inner sep=1.2pt,fill=black,draw,font=\sffamily\Large\bfseries}]

  \node[main node] (1) []{};
  \node[main node] (2) [below of=1] {};
  \node[main node] (3) [below of=2] {};

  \path[every node/.style={font=\sffamily\small}]

    (2) edge [] node[] {} (1)
    (3) edge [] node[] {} (2);

  \node[] (caption) [below of=3,yshift=0.4cm,draw=none,fill=none] {$\mathfrak C_1$};

\end{tikzpicture}
\hspace{0.6cm}
\begin{tikzpicture}[->,>=stealth',shorten >=1pt,auto,node distance=1cm,
  thick,main node/.style={circle,inner sep=1.2pt,fill=black,draw,font=\sffamily\Large\bfseries}]

  \node[main node] (1)[] {};
  \node[main node] (2) [right of=1] {};
  \node[main node] (3) [below of=1] {};
  \node[main node] (5) [below right of=3,xshift=-0.2cm] {};
  \node[main node] (4) [right of=3] {};

  \path[every node/.style={font=\sffamily\small}]

    (3) edge [] node[] {} (1)
    (3) edge [] node[] {} (2)
    (4) edge [] node[] {} (2)
    (4) edge [] node[] {} (1)
    (5) edge [] node[] {} (3)
    (5) edge [] node[] {} (4);

  \node[] (caption) [below of=5,yshift=0.4cm,draw=none,fill=none] {$\mathfrak C_2$};

\end{tikzpicture}
\par\end{center}
Now assume $G_2\neq\emptyset$.
Then $\mathfrak G$ is a rooted frame without clusters (otherwise it would subreduce to either $\mathfrak B_1$ or $\mathfrak B_2$), with at most $2$ strict successors from  each $G_1$ point (otherwise it would be subreducible to the trident $\mathfrak B_3$) and all $G_2$
points are only related to themselves (otherwise it would be reducible to $\mathfrak B_4$). 
\begin{center}

\begin{tikzpicture}[->,>=stealth',shorten >=1pt,auto,node distance=1cm,
  thick,main node/.style={circle,inner sep=1.2pt,fill=black,draw,font=\sffamily\Large\bfseries}]

  \node[main node] (1) []{};
  \node[main node] (2) [below of=1] {};
  \node[main node] (3) [below of=2] {};
  \node[main node] (4) [below of=3] {};

  \path[every node/.style={font=\sffamily\small}]

    (2) edge [] node[] {} (1)
    (3) edge [] node[] {} (2)
    (4) edge [] node[] {} (3);

  \node[] (caption) [below of=4,yshift=0.4cm,draw=none,fill=none] {$\mathfrak B_4$};

\end{tikzpicture}
\end{center}
We claim that, $G_1 \cup G_2$ remains connected after removing of the root of the frame, since otherwise it would be subreducible to the forbidden frame $\mathfrak B_5$. Indeed, suppose for sake of contrary that $G_1 \cup G_2$ is disconnected in a sense that two points $x,y\in G_1 \cup G_2$ are such that $x$ is not reachable from $y$ by the reflexive-transitive closure of the relation $R\cup R^{-1}$. Then $G_1 \cup G_2$ breaks down into finitely many connected components, at least one of which is of depth $1$ (since we assumed that $G_2\neq\emptyset$). Let $U$ be a component of depth $1$ and let $V$ be its complement inside of $G_1 \cup G_2$. Then we can map $G_1\cap U$ to the point $2$ of the
frame $\mathfrak B_5$,  $G_2\cap U$ to the point $1$, and all of the $V$ to the point $4$. Send the root $r$ to the point $3$. It is easy to check that the map thus defined is a p-morphism from $\mathfrak G$ to $\mathfrak B_5$. 

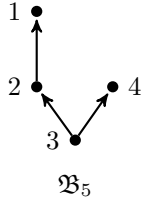
\begin{center}
\begin{tikzpicture}[->,>=stealth',shorten >=1pt,auto,node distance=1cm,
  thick,main node/.style={circle,inner sep=1.2pt,fill=black,draw,font=\sffamily\Large\bfseries}]

  \node[main node,label=left:{1}] (1)[] {};
  \node[main node,label=left:{2}] (2) [below of=1] {};
  \node[main node,label=left:{3}] (3) [below right of=2,xshift=-0.2cm] {};
  \node[main node,label=right:{4}] (4) [right of=2] {};

  \path[every node/.style={font=\sffamily\small}]

    (2) edge [] node[] {} (1)
    (3) edge [] node[] {} (2)
    (3) edge [] node[] {} (4);

  \node[] (caption) [below of=3,yshift=0.4cm,draw=none,fill=none] {$\mathfrak B_5$};

\end{tikzpicture}
\par\end{center}
This is a contradiction showing that the 'upper part' of $\mathfrak G$ is connected.
Consider $G_1 \cup G_2$ with a relation $R\cup R^{-1}$ as an undirected graph $Gr$.
Connectedness says that there exists a closed path that includes all
edges from $Gr$. Indeed, if we enumerate all the edges $e_1,\ldots,e_k$, then
there exists a path that connects one vertex of $e_i$ with another vertex of $e_{i+1}$. The concatenation of all such paths, with $e_i$ between, and ending with a path from a vertex of $e_k$ to a vertex of $e_1$ is a closed path that traverses all the edges. Since at least one such path exists, we may choose a shortest one, call it $l$. Let's write $l$ as
a word $x_0,\ldots, x_0$ starting from some point $x_0\in G_2$. It is clear that $x_i\overline{R}x_{i+1}$ or $x_{i+1}\overline{R}x_{i}$ for each $i<m$.  If $x_i\in G_1$ then we have $x_i\overline{R}x_{i-1}$ and $x_i\overline{R}x_{i+1}$. It may happen, however, that even though $x_i$ has two distinct successors in $G_2$ (recall that it cannot have more than two), still $x_{i-1}=x_{i+1}$. We wish to resolve such ``defects''. If we have $yxy$ somewhere in $l$, with $x\in G_1$, and $x\overline{R}z\neq y$, then we just replace it with $yxzxy$. This operation resolves current ``defect'' and does not add new ones, so after a finite number of steps we'll get a path $l\sp{\prime}$ which 
\begin{itemize}
\item[(a)]{traverses all edges in $Gr$;}  
\item[(b)]{ is free from the ``defects'', i.e. for any point $x\in G_1$, if $x$ has two distinct proper successors $y\neq z$, then whenever $x$ appears in $l\sp{\prime}$, it is preceded by one of the $y,z$ and followed by the other.}
\end{itemize}
Let $l\sp{\prime}$ be written as a word $x_0,x_1,\ldots,x_m,x_0$.
Define $S$ as set of all non-empty prefixes of $l\sp{\prime}$.

Now we are ready to define the desired $p$-morphic preimage $\mathfrak F$ of our
frame $\mathfrak G$ and the p-morphism. The intuition is that each path(point) in $\mathfrak F$ is the pre-image
of its last element.
Define $\mathfrak F=(W\sp{\prime},Q)$ as follows:  $W\sp{\prime} = \{r\} \cup S$ and $Q$ is the reflexive closure of the following.
\begin{itemize}
\item $r Q a$ for all $a\in S$.
\item $aQb$ if $last(a)~R~last(b)$ and  $b=ax$ or $a=bx$, with  $x\in G_1 \cup G_2$.
\end{itemize}

Here by $last(a)$ we mean the last element of the sequence $a$.

Let us define $g:W\sp{\prime}\to W$ as $g(x)=r$ if $x=r$ and $g(a)=last(a)$ otherwise.

We now show  that (1) $\mathfrak F$ is a crown frame and (2)
 $g$  is a p-morphism.
\begin{itemize}
\item[(1.)]{follows from the definition of $\mathfrak G$ and $\mathfrak F$. Indeed, let the even prefixes $x_0,x_1,\ldots,x_{2i}$ of $l\sp{\prime}$ be called $s_{2i}$ and the odd prefixes $x_0,x_1,\ldots,x_{2i+1}$ be called $s_{2i+1}$. It is clear from the way we constructed $l\sp{\prime}$ that $s_{2i}Q s_{2i+1}$ and $s_{2i}Qs_{2i-1}$. Moreover, since the states in $l\sp{\prime}$ alternate between $G_1$ and $G_2$, it is evident that the length of $l\sp{\prime}$ is odd (counting the last $x_0\in G_2$). Compare this with the presentation of the crown frame given in definition \ref{saw} to make sure that $\mathfrak F$ is isomorphic to a crown frame.}

\item[(2.)]{We verify that $g$ is a p-morphism. That $g$ is monotone is almost immediate from the definition.
To verify the back condition: for root it holds  trivially because the frames are transitive.
For the points $s_{2i}$ there is nothing to check as they are only related to themselves.
Let  $g(a)=x\in G_1$ and $x\overline{R}y$. Then $a=cx$ for some $l'$-prefix $c$.
 As $x\overline{R}y$, by construction of $l'$, $x$ and $y$ are always neighbors in $l'$, so
 there also exists a $l'$-prefix $d$ such that $g(dy)=y$ and either $c=dy$  or
 $d=cx$. By definition then $cxQdy$, as desired.}
 
 \end{itemize}

Consequently, the crown frame $\mathfrak F$ is reducible to $\mathfrak G$.
\end{proof}

We note in passing that the subframe $G_1\cup G_2$ considered in the last part of the above proof can be represented as a connected graph in yet another way. Namely, we can take the points of $G_2$ as vertices and the points in $G_1$ as edges between those points of $G_2$ that they relate to. This turns $G_1\cup G_2$ into a connected graph $G'$. In such a way, finding a coherent path like $l\sp{\prime}$ amounts to the so called Chinese Postman Problem for $G'$. Efficient algorithms are known for this problem, which one can employ to feasibly build a small suitable crown frame. Here we only concerned ourselves with the existence proof, since that is sufficient for our axiomatization task.

\begin{proposition}\label{SecondCompl}
The logic $\mathbf{PL}_2$ is axiomatized by the formula $\xi$ above $\mathbf{S4}$. In other words, $\mathbf{PL}_2=\Xi$.
\end{proposition}
\begin{proof}
By Theorem~\ref{t:Xi fmp} the logic $\Xi$ has the finite model property. Since any rooted finite frame for $\Xi$ is a reduction of a crown frame for $\mathbf{PL}_2$ by Lemma~\ref{l:crowns reduce} we have $\Xi\subseteq \mathbf{PL}_2$. By Lemma~\ref{l:crowns validate} we also have that each crown frame validates $\xi$, hence $\mathbf{PL}_2\subseteq\Xi$. It follows that $\mathbf{PL}_2=\Xi$.
\end{proof}

We also present a slightly more intuitive and concise axiomatization of $\mathbf{PL}_2$ by the following two formulas:
\[
\begin{array}{cl}
\textrm{(I)}& p \rightarrow \Box[\neg p \rightarrow \Box(p \rightarrow \Box p)]
\\
\textrm{(II)} & \Box[(r\wedge q ) \rightarrow \gamma ] \rightarrow [(r\wedge q) \rightarrow \di( \neg(r\wedge q) \wedge \di\Box p \wedge \di\Box\neg p)]
\end{array}
\]

Where $\gamma$ is the formula
\[
\di\Box(  p \wedge q)  \wedge  \di\Box(  \neg p \wedge q) \wedge \di\Box(  p \wedge \neg q).
\]
The formula (I) forbids frames $\mathfrak B_1$, $\mathfrak B_2$ and $\mathfrak B_4$, while (II) forbids $\mathfrak B_3$ and $\mathfrak B_5$. All crown frames validate both (I) and (II), thereby proving:

\begin{theorem}\label{ThirdCompl}
The logic $\mathbf{PL}_2$ is axiomatized by \emph{(I)} and \emph{(II)} over $\mathbf{S4}$. In other words, $\mathbf{PL}_2=\mathbf{S4}+\textsc{(I)}+\textsc{(II)}$.
\end{theorem}

Note that (I) carries an interesting dimensional meaning. Denote by $\delta A=\mathbb{C} A{-}A$ the \emph{external boundary} of $A$ (closure of $A$ minus $A$). Then a space $X$ validates (I) iff $\delta^3 A=\emptyset$ for all $A$. If $A$ is a polygon, then $\delta A$ is a polygon of strictly lower dimension. So over the polygonal plane $\delta^3 A=\emptyset$.

\section{Complexity of the satisfiability problem.}\label{complexity}
From the fmp and the finite axiomatization we conclude that our logic $\mathbf{PL}_2$ is
decidable. Moreover, we calculate the computational complexity of the satisfiability problem.

\begin{theorem}
The satisfiability problem of  our logic is $\PSPACE$ complete.
\end{theorem}

\begin{proof}
The encoding in \cite{wolt:spat00} shows that the following problem is $\PSPACE$-hard:
for $\phi$ and $\psi$ formulas, can $\phi$ be satisfied on a saw-model in which $\psi$ is true in every world? The saw-models in \cite{wolt:spat00} have the form $m_1e_1m_2e_2\ldots m_n$, where the states are the letters, and $R$ is the reflexive closure of the pairs $(m_i,e_{i+1})$ and $(m,e_{i-1})$.

Let $r,m,e$ be propositional variables. Let $C$ be the conjunction of these formulas:
\begin{eqnarray}
\  & \ & \mbox{$r,m,e$ are disjoint and one of them holds at each world.}\\
r &\rightarrow& \di m\\
m &\rightarrow& \di e.
\end{eqnarray}
In each model satisfying these formulas everywhere, $r$ is true only at the root, $e$ is true at all edge worlds, and the middle worlds make $m$ or $e$ true.

It is straightforward to show that the following are equivalent:
\begin{itemize}
\item $\phi$ is satisfied on a saw-model in which $\psi$ is true in every world;
\item $r\wedge \Box C \wedge \di((m\vee e)\wedge \phi) \wedge \Box((m\vee e)\rightarrow \psi)$ is satisfiable in our logic.
\end{itemize}
The upper bound algorithm uses a similar divide and conquer strategy for checking connectedness as used in \cite{wolt:spat00}.

Fix a formula $\theta$.

Let $\subf{\theta}$ be the smallest set containing all subformulas of $\theta$ and being closed under single negations.

A \emph{mosaic} is a structure $(W,R,l)$ with $W=\{r,m,e_0,e_1\}$, $R=\{(e,m),(m,e_0),(m,e_1)\}$ and $l:W\mapsto \Power{\subf{\theta}}$.

A mosaic is \emph{coherent} if it satisfies the following:
\begin{description}
\item[Bool] Each $l(w)$ is a maximal consistent subset of  $\subf{\theta}$.
\item[Box]  If $\di\phi\in \subf{\theta}$ and $vR^*w$ and $\phi \in l(w)$, then $\di\phi\in l(v)$.
\item[Middle-Di] If $\di\phi\in l(m)$, then $\phi\in l(e_0)$ or $\phi\in l(e_1)$.
\item[Edge-Di] If $\di\phi\in l(e_i)$, then $\phi\in l(e_i)$.
\end{description}
A set $M$ of mosaics is \emph{saturated} if it satisfies the following:
\begin{description}
\item[Root] All roots in all mosaics in $M$ are labeled by the same set of subformulas.
\item[Witness ]  If $(W_1,R_1,l_1)\in M$ and $\di\phi\in l(r_1)$, then there exist  a $(W_2,R_2,l_2)\in M$ and some $w\in W_2$ with $\phi\in l_2(w)$.
\item[Paths] If $m,m'\in M$, then there are $m=m_0,m_1,\ldots m_n=m'$ in $M$ which form a path. \\
We say that a string of mosaics $m_0,\ldots m_n$ form a path if $l_i(e_1^i)= l_{i+1}(e_0^{i+1})$.
\end{description}
From a model we can generate a saturated set of coherent mosaics, and conversely, we can build a model from such a set. Thus we obtain:

\begin{claim}
\label{thm:smallmodel}
A formula $\theta$ is satisfiable  if and only if there exists a saturated set of coherent mosaics with $\theta$ in the label of the root.
\end{claim}
As a corrollary we obtain that each formula can be satisfied in a model whose size is at most exponential in the length of the formula.

We now describe a $\PSPACE$ procedure which, given an input formula $\theta$, decides whether a   saturated set of coherent mosaics with $\theta$ in the label of the root exists.
The procedure is like a tableaux algorithm, but instead of working with sets of formulas, we work with mosaics.
\begin{description}
\item[Input] A formula $\theta$.
\item[Step root] Guess a coherent mosaic $m_0$ with $\theta$ in the label of the root $r$.
\item[Step witnesses] For each $\di\phi \in l(r)$, guess a coherent mosaic with the same root label as $l(r)$ and with $\phi$ in the label of one of the worlds. \\
We now have less than $|\theta|$ many mosaics, for which we can guess an order, say $m_0,m_1,\ldots m_n$. We can assume they will be glued together to a model in that order.
\item[Step path check] Now we must check that for each pair $m_i,m_{i+1}$,  
a path of mosaics from the first to the second element  exists. \\This is done with the    procedure $\checkpath(m,m')$, recursively defined as follows:
\\
$\checkpath(m,m')$ is true if either $l_m(e_1^m)= l_{m'}(e_0^{m'})$  (that is, they can be glued together) or there exists a coherent mosaic $m''$  whose root has the same label as the root of $m$ and is such that $\checkpath(m,m'')$ and $\checkpath(m'',m')$ are true.
\end{description}
The procedure outputs  that $\theta$ is satisfiable only if each step succeeds.
The algorithm is correct by Claim~\ref{thm:smallmodel}. It runs in non-deterministic polynomial space because the size of each mosaic is polynomial in the length of the input formula, in the witness step we have at most as many mosaics as the length of the input formula, and by always guessing a mosaic in the middle, the number of nested recursive calls of  the procedure $\checkpath$ is bounded by the length of the input formula, and thus can also be implemented in non-deterministic polynomial space. As non-deterministic $\PSPACE$ equals $\PSPACE$ by Savitch' Theorem, the procedure runs in $\PSPACE$.

\end{proof}

\section{Craig Interpolation}\label{interpolation}

Only very few extensions of $\mathbf{S4}$ have the Craig interpolation property \cite{maks:maks79}. Our logic is not among them, as the following counterexample shows.
Consider the following two formulas:

\[
\begin{array}{cl}
(A) & \Box(r \rightarrow \di(\neg r \wedge p \wedge \di\neg p)) \\
(C)& (r \wedge \di\Box s \wedge \di\Box\neg s) \rightarrow \di(\neg r \wedge \di\Box s \wedge \di\Box\neg s)
\end{array}
\]
 $A$ and $C$ have only the variable $r$ in common.
We claim that  (1)  $A \rightarrow C$ is valid, but (2) there is no interpolant for $A \rightarrow C$. That is, there is no formula $I$ written only in the variable $r$ such that $A\rightarrow I$ and $I \rightarrow C$ are both valid.

In the proof we use the notion of $\Sigma$-bisimulations, for $\Sigma$ a set of variables. These are bisimulations which only preserve the variables in $\Sigma$, not all variables in the language.

(1) Take any model $M$ satisfying $A$ and the antecedent of $C$. Then $M$ must be of depth 3,  $r$ is true only at the root, and  there are  endpoints making $s$ and $\neg s$ true. But then there must be a predecessor of an $s$ and a $\neg s$ end-point. As $r$ is only true at the root, this predecessor is the state asked for in the consequent of $C$.

(2) Take any model $M_1$ satisfying $A$ and $r$ at the root. Then $r$ is true only at the root. Thus $M_1$ $\{r\}$-bisimulates with the model $M_0$ consisting of a root and one succesor in which $r$ is only true at the root.
Let $M_2$ be the model consisting of one root and two successors which are end-points. In $M_2$  $r$ is true only at the root, and one of the end-points makes $s$ true. Then $M_2$ also $\{r\}$-bisimulates with $M_0$. Clearly $C$ is false at the root of $M_2$. Now assume to the contrary that $I$ is an interpolant. Then, because $A$ is true at the root of $M_1$, $I$ is true at the root of $M_1$, and thus by the $\{r\}$-bisimulation, $I$ is true at the root of $M_0$, and also at the root of $M_2$. But then $C$ must be true at the root of $M_2$, a contradiction.

\section{First order definability}\label{frstrddfnblt}
After all we prove that class of crown frames respect to which the modal logic $\mathbf{PL}_2$ is sound and complete, is not definable by means of first-order language.

\begin{theorem}
The class of \emph{crown frames} is not first-order definable.
\end{theorem}
\begin{proof}
Let us recall that crown frames have one important property. Namely, after removing the root the crown frame remains connected, and has finite-length path between any two possible worlds.

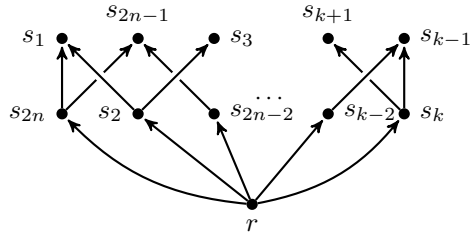
\begin{center}
\begin{tikzpicture} [->, >=stealth', shorten >=1pt,auto, node distance=1cm,
 thick,
 upa node/.style={circle,inner sep=1.2pt, fill=black, draw, font=\sffamily \Large \bfseries, label={left:$s_1$}},
 upb node/.style={circle,inner sep=1.2pt, fill=black, draw, font=\sffamily \Large \bfseries, label={$s_{2n-1}$}},
 upc node/.style={circle,inner sep=1.2pt, fill=black, draw, font=\sffamily \Large \bfseries, label={right:$s_3$}},
 upd node/.style={circle,inner sep=1.2pt, fill=black, draw, font=\sffamily \Large \bfseries, label={$s_{k+1}$}},
 upe node/.style={circle,inner sep=1.2pt, fill=black, draw, font=\sffamily \Large \bfseries, label={right:$s_{k-1}$}},
 da node/.style={circle,inner sep=1.2pt, fill=black, draw, font=\sffamily \Large \bfseries, label={left:$s_{2n}$}},
 db node/.style={circle,inner sep=1.2pt, fill=black, draw, font=\sffamily \Large \bfseries, label={left:$s_2$}},
 dc node/.style={circle,inner sep=1.2pt, fill=black, draw, font=\sffamily \Large \bfseries, label={right:$s_{2n-2}$}},
 dd node/.style={circle,inner sep=1.2pt, fill=black, draw, font=\sffamily \Large \bfseries, label={right:$s_{k-2}$}},
 de node/.style={circle,inner sep=1.2pt, fill=black, draw, font=\sffamily \Large \bfseries, label={right:$s_{k}$}},
 r node/.style={circle,inner sep=1.2pt, fill=black, draw, font=\sffamily \Large \bfseries, label={below:$r$}}]

  \node[upa node] (1) {};
  \node[upb node] (2) [right of=1] {};
  \node[upc node] (3) [right of=2] {};
  \node[upd node] (4) [right of=3,xshift=0.5cm] {};
  \node[upe node] (5) [right of=4] {};
  \node[da node] (6) [below of=1] {};
  \node[db node] (7) [below of=2] {};
  \node[dc node] (8) [below of=3] {};
  \node[dd node] (9) [below of=4] {};
  \node[de node] (10) [below of=5] {};
  \node[r node] (r) [below of=8,xshift=0.5cm, yshift=-0.2cm] {};

  \path[every node/.style={font=\sffamily\small}]

    (6) edge [] node[] {} (1)
        edge [] node[] {} (2)
    (8) edge [] node[] {} (2)
        edge [draw=none] node[] {$\cdots$} (9)
    (7) edge [-,draw=white, line width=4pt] node[] {} (1)
        edge [] node[] {} (1)
       edge [-,draw=white, line width=4pt] node[] {} (3)
        edge [] node[] {} (3)
    (10)edge [] node[] {} (4)
        edge [] node[] {} (5)
    (9) edge [-,draw=white, line width=4pt] node[] {} (5)
         edge [] node[] {} (5)
    (r) edge [] node[] {} (7)
        edge [] node[] {} (8)
        edge [] node[] {} (9);

       \draw[->] (r) to [bend right=20] node[] {} (10);
       \draw[->] (r) to [bend left=20] node[] {} (6);

\end{tikzpicture}
\end{center}
Note, that for any natural number $n\in \mathbb{N}$, there is exactly one crown frame $\SAW_n$ with  $2n+1$ possible worlds.  Hence let us number a family of all crown frames by means of $\mathbb{N}$, the set of natural numbers. Let $s^{n}_{i}$ denotes a point $s_i$ in $n$-th crown frame $\SAW_n$. Let us consider $f(n)=s_{1}^{n}$ and $g(n)=s_{n+1}^{n}$, two members of $\prod_{n\in \mathbb{N}}\SAW_n$, and a non-principal, that is free ultrafilter $U$ on $\mathbb{N}$. We claim, that there is no finite path from equivalence class of $f$ to equivalence class of $g$ in the ultraproduct $\prod_{U}\SAW_n$. For sake of contrary suppose, that there is a finite $k$-tuple $h_{1}, h_2, \cdots , h_{k}$ such that $[f]=[h_{1}], [h_2], \cdots , [h_{k}]=[g]$ is a  path from equivalence class of $f$ to equivalence class of $g$. Then there is an element of ultrafilfer $A$ (since path is finite) such that for every $1\leq i,j\leq k$, $A\subseteq \{n| ~h_{i}(n)Q_{n} h_{j}(n)\}$  or $A\subseteq \{n| ~h_{j}(n)Q_{n} h_{i}(n)\}$. Note that since an ultrafilter  $U$ is non-principal then it doesn't contain any finite subset as its own element. Since length of path from $s_{1}^{n}$ to $s_{n+1}^{n}$ increases with $n$ we get contradiction. Then the class of crown ~frames isn't closed under ultraproducts and hence by Proposition \ref{fod} the class is not first order definable.
\end{proof} 
\newpage
\section*{References}
\bibliographystyle{elsarticle-harv}
\bibliography{pl2}
\end{document}